\documentclass[a4paper,12pt]{article}

\usepackage{citesort,geometry,enumerate,amsmath,amssymb,latexsym}
\newtheorem{example}{Example}[section]
\newtheorem{Def}[example]{Definition}
\newtheorem{Prop}[example]{Proposition}
\newtheorem{Theo}[example]{Theorem}

\newtheorem{Cor}[example]{Corollary}

\newtheorem{Exam}[example]{Example}
\newenvironment{Prf}{{\bf Proof:}}{\hfill $\Box$

\mbox{}}

\def\wtilde{\widetilde}

\begin{document}

\title{\Large\bf Extendibility , monodromy and  local triviality
 for topological groupoids 
}
\small{ \author  { Osman Mucuk \\ University of Erciyes\\
Faculty of Science and Art \\Department of Mathematics\\Kayseri/Turkey\\
mucuk@erciyes.edu.tr \\ \and \.{I}lhan \.{I}\c{c}en    \\  University of  \.{I}n\"{o}n\"{u} \\
Faculty of Science and Art\\ Department of Mathematics
\\ Malatya/ Turkey \\ iicen@inonu.edu.tr }}
\maketitle
\begin{abstract}
\footnotesize{A groupoid is a small category in which each morphism has an
inverse. A topological groupoid is a groupoid in which both sets of objects and
morphisms have topologies such that all groupoid structure maps are continuous. 
The notion of monodromy groupoid of a topological groupoid  generalises those of 
fundamental groupoid and universal covering.  It was  earlier proved that  the monodromy 
of a locally sectionable topological groupoid has a topological groupoid structure satisfying 
some properties. In this paper a similar problem is studied for compatible locally trivial 
topological groupoids.}
\end{abstract}

\section*{Introduction}
Let $G$ be a topological groupoid and $W$ an open subset of $G$ containing all the identities. 
Then the monodromy groupoid $M(G,W)$ of $G$ is defined as in Definition \ref{mon}, which is due
 to Pradines \cite{Pr}. In \cite{Mu} 
(see also \cite{BM1} and
\cite{BM2} for Lie groupoid case ) in the case where $G$ is locally sectionable it was proved 
that the groupoid $M(G,W)$  may be given
a topology making it a topological groupoid  such that each star $M(G,W)_x$ is a universal covering 
of $G_x$ and $M(G,W)$ has a universal property on the globalisation  of continuous pregroupoid
morphisms to  topological groupoid morphisms. 

If $X$ is a topological space  then $G=X\times X$ is a topological groupoid with the groupoid
multiplication $(x,y)(y,z)=(x,z)$.  If  $X$ is a path connected topological space which has a 
universal covering then the monodromy groupoid of the topological groupoid  $G=X\times X$ for
a suitable open neighbourhood $W$ is the  fundamental groupoid $\pi_1(X)$.
 
If $G$ is a topological group, which can be thought as a topological groupoid with only one object, 
such that it has a universal covering, then for  a suitable open subset $W$ of  $G$ the monodromy 
groupoid $M(G,W)$ is  just the universal covering of $G$. 

Thus the concept of the  monodromy groupoid  generalises the notions of the fundamental groupoid as 
topological groupoid and the universal covering.

In \cite{Ma} the monodromy groupoid denoted by $\Pi G$ of a topological groupoid $G$ in which each
 star $G_x$ has a universal covering is constructed directly from the universal coverings of 
 ${G_x}'s$, and then $\Pi G$ is given a topology under some conditions   studying in terms of 
 principal bundles.

In this paper , we define compatible locally trivial groupoid and as similar to  the fundamental 
groupoid case studied  in \cite{BrDa} prove that the monodromy groupoid of a compatible locally 
trivial groupoid can be topologised with an appropriate  topology.

\section{Groupoids and topological groupoids}
 A {\em groupoid} $G$ on $O_G$ is a small category in which each morphism is
 an isomorphism. Thus $G$ has a set of morphisms, which we call just
 {\em elements} of $G$, a set $O_G$ of {\em objects} together with functions
 $\alpha, \beta\colon G\rightarrow O_G$, $\epsilon \colon O_G \rightarrow G$
 such that $\alpha\epsilon=\beta\epsilon=1_{O_G}$, the identity map. The functions
 $\alpha$, $\beta$ are called {\em initial} and {\em final} point maps respectively.
 If $a,b\in G$ and $\beta(a)=\alpha(b)$, then the {\em product} or {\em composite}
 $ab$ exists such that $\alpha(ab)=\alpha(a)$ and $\beta(ab)=\beta(b)$. Further,
 this composite is associative, for $x\in O_G$ the element $\epsilon (x)$ denoted
 by $1_x$ acts as the identity, and each element $a$ has an inverse $a^{-1}$ such
 that $\alpha(a^{-1})=\beta(a)$, $\beta(a^{-1})=\alpha(a)$, $aa^{-1}=(\epsilon \alpha)(a)$,
 $a^{-1}a=(\epsilon \beta)(a)$. The map $G\rightarrow G$, $a\mapsto a^{-1}$
 is called the {\em inversion}.

  In a groupoid $G$ for $x,y\in O_G$ we write $G(x,y)$ for the set of all
 morphisms with initial point $x$ and final point $y$. For
 $x \in O_G$ we denote the star \{$a\in G\colon\alpha(a)=x$ \} of $x$ by
 $G_x$ and the {\em costar} \{$a\in G\colon\beta(a)=x$ \} of $x$  by $G^x$.
In $G$ the set $O_G$ is mapped bijectively to the set of identities by
$\epsilon \colon O_G \rightarrow G$. So we sometimes write $O_G$ for the set of identities.
Let $G$ be a groupoid and $W$ a subset of $G$ such that $W\subseteq O_G$.
We say $G$ is {\em generated} by $W$ if each element of $G$ is written as a union
of some elements in $W$.

 Let $G$ be a groupoid. A {\em subgroupoid} of $G$ is a pair of subsets
 $H\subseteq G$ and $O_H\subseteq O_G$   such that $\alpha(H)\subseteq O_H$, $\beta(H)\subseteq O_H$,
 $1_x\in H$ for each $x\in O_H$ and $H$ is closed under the partial
 multiplication and inversion in $G$.

 A {\em morphism} of groupoids $H$ and $G$ is a functor, that is, it consists
 of a pair of functions $f\colon H\rightarrow G$ and $O_f\colon O_H \rightarrow O_G$
 preserving all the structures.
\begin{Def}{\em ( ~\cite{Ma}) {\em A topological groupoid} is a groupoid G on $O_G$ together with
topologies on G and $O_G$ such that the maps which define the groupoid
structure are continuous, namely the initial and final point maps
$\alpha,\beta\colon G\rightarrow O_G$, the object inclusion map
$\epsilon\colon O_G\rightarrow G$, $x\mapsto \epsilon (x)$, the inversion
$G\rightarrow G$, $a\mapsto a^{-1}$ and the partial multiplication
$ G_{\alpha}\times_{\beta}G\rightarrow G, (a,b)\mapsto ab$, where the pullback 
\[G_{\beta}\times_{\alpha} G=\{ (a,b)\in G\times G \colon \beta(a)=\alpha(b)\}\]
has the subspace topology from $G\times G$.}
\end{Def}
A {\em morphism} of topological groupoids $f\colon H\rightarrow  G$ is a 
morphism of groupoids in which both maps  $f\colon H\rightarrow G$
and  $O_f\colon O_H \rightarrow O_G$ are continuous.

Note that in this definition the partial multiplication $G_{\beta}\times_{\alpha}G\rightarrow G, 
(a,b)\mapsto ab$ and the inversion map  $G\rightarrow G, a\mapsto a^{-1}$ are continuous if and only if
 the map $\delta\colon G\times_{\alpha}G\rightarrow G, (a,b)\mapsto a^{-1}b$, called
the groupoid difference map, is continuous, where the pullbak
\[G\times_{\alpha} G=\{(a,b)\in G\times G \colon \alpha(a)=\alpha(b)\}\]
has the subspace topology from $G\times G$. Again if one of the  maps $\alpha,\beta$ and the
 inversion are continuous, then the other map is continuous.

Let $X$ be a topological space. Then  $G=X\times X$ is a topological
groupoid on $X$. In which each pair $(x,y)$ is a morphism from $x$ to $y$
and the groupoid composite is defined by $(x,y)(y,z)=(x,z)$. The inverse
of $(x,y)$ is $(y,x)$ and the identity at $1_x$ is the pair $(x,x)$. 

Note that in a topological groupoid $G$ for each $a\in G(x,y)$   right
translation $R_a\colon G^x\rightarrow G^y, b\mapsto ba$ and left translation
$L_a\colon G_y\rightarrow G_x$, $b\mapsto ab$ are homeomorphisms.

A groupoid $G$ in which each star $G_x$ has a topology such that  for 
$a\in G(x,y)$ the  right translation $R_a\colon G^x\rightarrow G^y, b\mapsto ba$ (and hence the
left translation
$L_a\colon G_y\rightarrow G_x, b\mapsto ab$) is homeomorphisms, is called 
{\em star topological groupoid}

Let $G$ be a groupoid and $W$ a subset of $G$ such that $O_G\subseteq W$. So the
set $W$ inherits a pregroupoid  structure from the groupoid $G$. That is the source and target
maps $\alpha$ and $\beta$ restrict to $W$ and if $a,b\in W$ and $\beta a=\alpha b$, then the
composition $ab$ of $a$ and $b$ in $G$ may or not belong to $W$. We follow the method of Brown
and Mucuk \cite{BM1}, which generalises work for groups in Douady and Lazard \cite{DL}. 

There is a standard construction $M(G,W)$ associating to the pregroupoid $W$ a morphism
$\tilde\imath\colon W\rightarrow M(G,W)$ to a groupoid $M(G,W)$ and which is universal for 
pregroupoid morphism to a groupoid. First we form th free groupoid $F(W)$ on the graph $W$
and denote the inclusion $W\rightarrow F(W)$ by $a\mapsto [a]$. Let $N$ be the normal
subgroupoid (Higgins\cite{Hi}, Brown\cite{Br}) of $F(W)$ generated by the elements 
$[a][b][ab]^{-1}$ for all $a,b\in W$ such that $ab$ is defined and belongs to $W$. Then 
$M(G,W)$ is defined to be the quotient groupoid $F(W)/N$. The composition $W\rightarrow
F(W)\rightarrow M(G,W)$ is written $\wtilde\imath$, and is the required universal morphism.
It is followed that $\wtilde\imath$ is injective.

The map $\wtilde\imath$ has a universal property that if $f\colon W\rightarrow H$ is a 
pregroupoid morphism  then there is a morphism of groupoids $f'\colon M(G,W)\rightarrow H$
such that $f= f'\wtilde\imath$. In particular the inclusion map $i\colon W\rightarrow G$
 globalises to a morphism  of groupoids $p\colon M(G,W)\rightarrow G$ called canonical
 morphism. 
 
 We give this construction as a definition.

\begin{Def}\label{mon}
{\em Let $G$ be a topological groupoid and  $W$ an open subset of $G$ such that
 $O_G\subseteq  W$.  Let $F(W)$ be the free groupoid on $W$ and let $N$ the normal subgroupoid
of $F(W)$ generated by the elements in the form $[a][b][ab]^{-1}$  for $a,b\in W$ such that
$ab$ is defined and $ab\in W$.  The quotient groupoid $F(W)/N$ is called 
{\em monodromy groupoid}  of $G$ for $W$and denoted by $M(G,W)$. }
\end{Def}
\section{Monodromy groupoids of locally sectionable topological groupoids}
In this section we recall some results from \cite{BM1}.

 Let  $G$  be a topological  groupoid such that each star $G_x$ has a universal covering.   
 The  groupoid $\Pi G$  as the union of the universal coverings of ${G_x}'s$ is  defined  
as follows.  As a set,  $\Pi G$  is the
union of the stars $(\pi_{1}G_x)_{1_x}$.  The  object  set of $\Pi
G$  is the same as that of  $G $.  The  function   $\alpha \colon
{\Pi G}\rightarrow X$  maps all of $(\pi_{1}G_x)_{1_x}$ to   $x $,
while   $\beta \colon  {\Pi G}\rightarrow X$   is  on
$(\pi_{1}G_x)_{1_x}$    the composition of the two target maps 
\[(\pi_1G_x)_{1_x}\rightarrow G\rightarrow X.\]

 As explained in Mackenzie~\cite{Ma}, p.67,
there is  a multiplication on   $\Pi G$ given by `concatenation',
i.e.
\[    [a] \circ [b] = [a+a(1)b] ,\]                  where the
$+$  inside the bracket denotes the usual composition of  paths.
 Here  $a $ is assumed to be a  path  in   $G_x$   from  $1_x$
to   $a(1) $,   where $\beta (a(1)) = y $, say, so that $b$  is
a path in  $G_y $, and for  each  $t\in [0,1] $, the product
$a(1)b(t)$  is defined in  $G $, yielding a  path $b(a(1))$
from  $a(1)$  to  $a(1)b(1) $.  It is straightforward to prove
that in this  way   $\Pi G$   becomes a groupoid, and that the
final map of paths induces a  morphism  of groupoids $p \colon
{\Pi G}\rightarrow G $.    

       Let  $X$  be a topological space admitting a simply
connected cover.    A subset  $U$  of  $X$  is called {\em
liftable} if  $U$  is open, path-connected and the inclusion
$U\rightarrow X$ maps each fundamental group of $U$  trivially.
If  $U$  is liftable, and  $q\colon Y\rightarrow X$  is a
covering map,  then  for  any   $y\in Y $  and  $x\in U$  such
that  $qy = x $, there is a unique map   $\hat{\imath}\colon U
\rightarrow Y$   such  that  $\hat{\imath} x = y$  and
$q\hat{\imath}$  is the  inclusion   $U\rightarrow X $.   This
explains  the  term {\it liftable}.    

\begin{Theo}\label{startop}   Suppose that  G  is a star connected
star  topological    groupoid  and  W  is an open neighbourhood of  $O_G$
satisfying the condition  

$(\star)$  W   is  star  path-connected  and  $W^2$    is
contained  in  a   star  path-connected neighbourhood  V  of
$O_G$  such that for all  $x\in O_G $,  $V_x$  is  liftable.

Then there is an
isomorphism  of star topological  groupoids  $M(G,W)\rightarrow
\Pi G $, and hence the morphism  $p\colon M(G,W)\rightarrow G$  is a
star universal covering map. 
\end{Theo}

The following definition is due to Ehresmann~\cite{Eh}.

\begin{Def} {\em Let  $G$  be a groupoid and let  $X = O_{G}$
be a topological space. An {\em admissible local section} of
$G$ is a function $s : U \rightarrow G$  from an open set in $X$
such that

1. $\alpha s(x) = x$ for all $x \in U$;

2.  $\beta s(U)$ is open in $X$, and

3.  $\beta s$ maps $U$  homeomorphically to $\beta s(U)$. }
\label{locsecdef}
 \end{Def}
Let $W$ be a subset of $G$ such that $X\subseteq W$ and let $W$ have the structure of a
topological space. We say that
$(\alpha ,\beta ,W)$  is {\em  locally sectionable} if for each
$w \in W$ there is an admissible  local section $s : U
\rightarrow G$ of $G$ such that (i) $s\alpha (w) = w$,  (ii)
$s(U) \subseteq W$ and (iii) $s$ is smooth as a function from
$U$  to $W$. Such an $s$ is called a smooth {\it admissible
local section.}\par

The following definition is due to Pradines~\cite{Pr} under the
name  ``{\it morceau de groupoide diff\'erentiables}''. 
\begin{Def} {\em A {\it
locally topological   groupoid} is a  pair $(G,W)$ consisting of a
groupoid $G$ and a topological space smooth  $W$ such that:

\noindent $G_1)$  \ \  $O_G \subseteq W \subseteq G $;

\noindent $G_2)$ \ \ $W = W^{-1} $;

\noindent $G_3)$ \ \  the set $W(\delta ) = (W \times _{\alpha}
W) \cap \delta ^{-1}(W)$  is open in $W \times _{\alpha} W$ and
the restriction of $\delta $ to  $W(\delta )$ is continuous;

\noindent $G_4)$ \  \ the restrictions to $W$ of the source and
target maps $\alpha $  and $\beta $ are continuous  and the triple
$(\alpha ,\beta ,W)$ is  locally sectionable;

\noindent $G_5)$ \ \ $ W$ generates $G$ as a groupoid.}
\end{Def}

Note that, in this definition, $G$ is a groupoid but does not
need  to have a topology. The locally Lie groupoid $(G,W)$ is
said  to be {\it extendible} if there can be found a topology on
$G$ making it  a topological groupoid and for which $W$ is an open
subset.  In general $(G,W)$ is not extendible, but there is a holonomy
groupoid $Hol(G,W)$ and a morphism $\phi\colon (G,W)\rightarrow G$ of groupoids such that
$Hol(G,W)$ admits the structure of topological groupoid. The construction is given
in detail in \cite{AoB}. For an example of locally topological 
groupoid which is not extendible see \cite{AoB}.

\begin{Theo}\label{glob} (Globalisability theorem)
 Let $(G,W)$  be a locally topological  groupoid. Then there is a
topological   groupoid $H$,  a morphism $\phi : H \rightarrow
G$   of groupoids and an embedding $i : W \rightarrow H$
 of $W$  to an open neighbourhood of $O_{H}$  such
that the following  conditions are satisfied.

\noindent i) $\phi $  is the identity on objects, $\phi i =
id_{W} ,  \phi ^{-1}(W)$ is open in $H$,and the
restriction  $\phi _{W} : \phi^{-1}(W)  \rightarrow W$  of
$\phi $  is continuous;

\noindent  ii) if A is a topological   groupoid and $\xi :  A
\rightarrow G$  is a morphism of groupoids such that:

a) $\xi $ is the identity on objects;

b) the restriction $\xi _{W} : \xi ^{-1}(W) \rightarrow W$ 
of  $\xi $  is smooth and $\xi ^{-1}(W)$  is open in A
and generates A;

$c)$  the triple $(\alpha _{A} , \beta _{A} , A)$  is
locally  sectionable,

\noindent  then there is a unique morphism $\xi ^\prime : A
\rightarrow H$   of topological  groupoids such that $\phi \xi
^\prime = \xi $   and $\xi ^\prime a = i\xi a$  for $a
\in \xi ^{-1}(W)$. 
\end{Theo}

The groupoid $H$ is called the {\it holonomy groupoid}
$Hol(G,W)$ of the  locally topological groupoid $(G,W)$.

Let $G$ be a topological groupoid on $X$. Then $G$ is called {\it
locally trivial} if for all $x\in X$ there is an open set $U$
containing $x$ and a  section $s\colon U\rightarrow G_x$ of $\beta$.
Thus $\beta s=1_{U}$ and  for each $y\in U$, $\alpha(s(y))=x$,
i.e. $s(y): x \to y$ in $G$. 

The following result, whose proof is due to Pradines, is given in \cite{AoB}. 
In the proof of this result the sectionable condition  is used.
\begin{Prop}{\bf( 6.1 in \cite{AoB})}
A locally trivial locally topological groupoid is extendible.
\end{Prop}

As a corollary of  of Theorem \ref{glob} the following result is 
 obtained. See \cite{BrIc} for a nice application of this result 
 in the holonomy groupoids of local subgroupoids.
 
\begin{Cor}\label{locexcor} {\it Let} $G$ {\it be} a topological  {\it
groupoid  and let} $p : M \rightarrow G$ {\it be a morphism of
groupoids such that} $p : O_{M} \rightarrow O_{G}$ {\it  is the
identity. Let} $W$ {\it be an open subset of} $G$ {\it such that}

$a) \ \  O_{G} \subseteq W$;

$b) \ \ W = W^{-1}$;

$c) \ \  W$ {\it generates} $G${\it ;}

$d) \ \ ( \alpha _{W}, \beta _{W}, W )$ {\it is  locally
sectionable;}

{\flushleft\it and suppose that} $\tilde{\imath} : W
\rightarrow M$ {\it is given  such that} $p \tilde{\imath} = i :
W \rightarrow G$ {\it is the inclusion and}  $W^\prime =
\tilde{\imath}(W)$ {\it generates} $M ${\it .}

{\it Then} $M$ {\it admits a unique structure of topological  groupoid
such that} $W^\prime $ {\it is an open subset and}  $p : M
\rightarrow G$ {\it is a morphism of topological groupoids  mapping}
$W^\prime $  {\it homeomorphically to} $W ${\it .}  
\end{Cor}

These imply
Theorem \ref{stmonp}  which is called the {\em Strong Monodromy
Principle}, namely the  globalisation of   smooth local
morphisms to smooth morphisms on the  monodromy groupoid.

\begin{Theo} Let  G  be a locally sectionable topological   groupoid and
let W  be an open subset  of   G   containing   $O_G$,  such
that   $W = W^{-1} $,  and   W  generates   G.   Then  the
monodromy  groupoid  $M = M(G,W)$   admits  the structure of topological 
groupoid such that  $\tilde{\imath}(W)$  is an open  subspace
of  M  and any  continuous   pregroupoid morphism on  W  globalises to a
continuous  morphism on  M.  \label{monglob} \end{Theo}

\begin{Theo} \label{moniso} Suppose further to the assumptions of
Theorem \ref{monglob} that  G  is star path-connected, that
each of its stars has  a  simply  connected  covering,  and that
$W^2$   is contained in an open neighbourhood  V  of  $O_G$
such  that each star  $V_x$  in  $G_x$  is liftable.  Then the
projection   $p \colon M(G,W) \rightarrow G$  is the universal
covering map on each star, and so   M(G,W)   is  isomorphic to
the  star universal cover $\Pi G$ of  G.   \end{Theo}

\begin{Theo}\label{stmonp} {\em (Strong monodromy principle)}
Let  $G$ be a locally sectionable star path-conn-ected Lie
groupoid and let $W$ be a neighbourhood of  $O_G$ in  $G$ such
that W satisfies the condition: \newline ($\star$) W is
path-connected and   $W^2$ is contained in a star path-connected
neighbourhood V of $O_G$ such that for all $x\in O_G$, $V_x$  is
liftable. \newline \noindent  Let $f: W \rightarrow H$ be a
continuous pregroupoid morphism from W to the topological  groupoid H. Then f
determines uniquely a morphism $f' :\Pi G \rightarrow H$ of Lie
groupoids such that $f'j' = f$.  \hfill $\Box$ \end{Theo}


\section{Extendibility of compatible locally trivial  groupoids}

The concept of locally triviality  is due to Ehresmann~\cite{Eh}. It is also used  by Mackenzie in 
\cite{Ma}. We adapt this concept as follows.

\begin{Def}\label{comlc}{\em
Let $X$ be a topological space and  let $G$ be a groupoid on $X$. Let $\mathcal U=\{U_i\colon i\in I\}$
be an open cover of $X$, which is also a base for $X$ such that if $x\in U_i$ then there is a local
section  $s_{x,i}\colon U_i\rightarrow G_x$ of $\beta$ such that $s_{x,i}(x)=x$. Thus for $y\in U_{x,i}$, 
$s_{x,i}(y)\colon x\rightarrow y$ in $G$. Such a map $s_{x,i}\colon U_i\rightarrow G_x$ is called 
 {\em local section} about $x\in O_G$.  $G$ is called {\em compatible locally trivial} if 
 the following compatible condition is satisfied:\\
{\bf Comp:} For $x\in X$ if $s_{x,i}\colon U_i\rightarrow G_x$ and 
$s_{x,j}\colon U_j\rightarrow G_x$ are two local sections about $x$ then there
 is an open neighbourhood $V_{ij}$ of $x$ in $\mathcal U$ such that $V_{ij}\subseteq U_i\cap U_j$ and
$s_{x,i}|V_{ij}=s_{x,j}|V_{ij}$.}
\end{Def}
 If $s_{x,i}\colon U_{i}\rightarrow G_x$ is a local section about $x$, then we write
$\wtilde{U}_{x,i}$ for the image $s_{x,i}(U_{i})$.
We now prove the following theorem  on the extendibility of compatible locally trivial 
groupoid to a topological groupoid. This result is similar to the fundamental groupoid case 
studied earlier in \cite{BrDa}.

\begin{Theo}\label{comploctr} Let G be a  compatible locally trivial groupoid.  Let $a\in G(x,y)$. Then 
the sets  $({\wtilde{U}_{x,i}})^{-1}a(\wtilde{U}_{y,j})$ for all $x\in U_{i}$ and 
$y\in U_{j}$ form a set of basic neighbourhoods for a topology on $G$ such that $G$ is a topological 
groupoid with this topology.
\end{Theo}
\begin{Prf}
By the compatibility condition it is obvious that these sets form a topology on $G$. 
Because if  $({\wtilde{U}_{x,i}})^{-1}a(\wtilde{U}_{y,j})$ and 
$({\wtilde{U}_{x,i'}})^{-1}a(\wtilde{U}_{y,j'})$ are basic open neighbourhood of $a$, then
$x\in U_i\cap U_{i'}$ and $y\in U_j\cap U_{j'}$. By the compatible conditions there is a basic open 
neighbourhood $V_{ii'}$ of $x$ in $\mathcal U$ such that
$V_{ii'}\subseteq U_i\cap U_{i'}$ and $s_{x,i}|V_{ii'}=s_{x,i'}|V_{ii'}$ . 
Similarly here is a basic open neighbourhood $V_{jj'}$ of $y$ such that
$V_{jj'}\subseteq U_j\cap U_{j'}$ and $s_{x,j}|V_{ij}=s_{x,j'}|V_{jj'}$. Hence 
$({\wtilde{V}_{x,ii'}})^{-1}a(\wtilde{V}_{y,jj'})$ is a basic open neighbourhood of $a$
such that
\[({\wtilde{V}_{x,ii'}})^{-1}a(\wtilde{V}_{y,jj'})\subseteq
 ({\wtilde{U}_{x,i}})^{-1}a(\wtilde{U}_{y,j})\cap ({\wtilde{U}_{x,i'}})^{-1}a(\wtilde{U}_{y,j'})\]

  We now prove that $G$ is a topological groupoid with this topology. To prove that 
the groupoid difference map
\[\delta\colon G\times_{\alpha} G\rightarrow G,(a,b)\mapsto a^{-1}b \]
 is continuous let $\delta(a,b)=a^{-1}b$, where $a\in G(x,y)$ and $b\in G(x,z)$. Let
 $({\wtilde{U}_{y,j}})^{-1}a^{-1}b({\wtilde{U}_{z,k}})$ be a basic neighbourhood of $a^{-1}b$. 
 Then $({\wtilde{U}_{x,i}})^{-1}a({\wtilde{U}_{y,j}}){\times_{\alpha}} 
 ({\wtilde{U}_{x,i}})^{-1}b({\wtilde{U}_{z,k}})$
is an open neighbourhood of $(a,b)$ and 
\[\delta(({\wtilde{U}_{x,i}})^{-1}a({\wtilde{U}_{y,j}}){\times_{\alpha}} 
 ({\wtilde{U}_{x,i}})^{-1}b({\wtilde{U}_{z,k}}))=
 ({\wtilde{U}_{y,j}})^{-1}a^{-1}b({\wtilde{U}_{z,k}}).\]
So $\delta$ is continuous.

To prove the continuity of the target point map $\beta\colon G\rightarrow X$ let 
$a\in G$ with $\alpha(a)=x$ and $\beta(a)=y$ and let $U_j$ be a basic open neighbourhood of $y$ in 
$\mathcal U$. Then 
$({\wtilde{U}_{x,i}})^{-1}a({\wtilde{U}_{y,j}})$
is a basic open neighbourhood of $a$ and $\beta(({\wtilde{U}_{x,i}})^{-1}a({\wtilde{U}_{y,j}}))\subseteq
U_j$. Hence $\beta$ is continuous. Further since $\delta$ and $\beta$ are both continuous so also is the
source point map $\alpha\colon G\rightarrow X$. 

Finally for the continuity of the identity point map
$\epsilon\colon X\rightarrow G$ let $x\in X$ and let $({\wtilde{U}_{x,i}})^{-1}1_x({\wtilde{U}_{x,j}})$
be a basic open neighbourhood of $1_x$. Then $x\in U_i\cap U_j$ and by the compatibility condition there
is an open neighbourhood of $x$ in $\mathcal U$ such that $V_{ij}\subseteq U_i\cap U_j$ and 
$s_{x,i}|V_{ij}=s_{x,j}|V_{ij}$. 
So $\epsilon (V_{i,j})\subseteq ({\wtilde{U}_{x,i}})^{-1}1_x({\wtilde{U}_{x,j}})$. Hence $\epsilon$ is
continuous.

\end{Prf}
In particular in this result if we take $G$ to be the fundamental groupoid $\pi_1X$ on a topological
space $X$ which is locally nice then we obtain a result given in \cite{BrDa} stated as follows:

\begin{Exam} {\em Let $X$ be a locally path connected and  semilocally 1-connected topological
space. Then there is an open cover $\mathcal U$ of $X$ as in Definition \ref{comlc} consisting of all 
open, path connected subsets $U$ of $X$ such that there is only one homotopy class of the paths in $U$ 
between same points. If $x\in U_i$ for $U_i\in \mathcal U$, then the local section 
$s_{x,i}\colon U_i\rightarrow \pi_1 X$ is defined  by choosing for each 
$y\in U_i$  a path in $U_i$ from $x$ to $y$ and taking $s_{x,i}(y)$ to be the homotopy class of this path. 
Since the path class of the paths in $U_i$ between same points is unique, the local section
$s_{x,i}$ is well defined. Further since $X$ is locally path connected the  compatibility condition
is satisfied. Thus  $\pi_1X$ is a compatible locally trivial groupoid.} 
\end{Exam}

Let $G$ be a topological groupoid and $W$ an open neighbourhood of $O_G$ such that 
$O_G\subseteq W$. Then we have the monodromy groupoid $M(G,W)$ defined as in
Definition\ref{mon}. Note that by the construction of $M(G,W)$ we have a pregroupoid morphism
$\wtilde\imath\colon W\rightarrow M(G,W)$ which is an  inclusion map. Let $\wtilde W$ be the image 
of $W$ under this inclusion and let 
$\wtilde W$ have the topology such that $\wtilde \imath\colon W\rightarrow \wtilde W$ is a homeomorphism. Then we have 
 following result.
\begin{Prop}\label{loctrmon}
Let $G$ be a topological groupoid  whose underlying groupoid  is compatible locally trivial and let $W$ be an 
open subset of $G$ such that $O_G\subseteq W$.   Suppose that each local section 
$s_{x,i}\colon U_{x,i}\rightarrow G_x$ has image in $W$, i.e. $s_{x,i}(U_{x,i})\subseteq W$. Then  the 
groupoid   $M(G,W)$ is compatible locally trivial.
\end{Prop}
\begin{Prf} The proof is immediate since $W$ is isomorphic to $\widetilde{W}$ and the inclusion
$\wtilde\imath\colon W\rightarrow \wtilde{W}$ is identity on objects.
\end{Prf}
As a corollary of this result we state the following.
\begin{Cor} Let G be a topological groupoid whose groupoid is compatible locally trivial and let $W$ be 
an open subset of $G$ including all the identities. Suppose that each local section 
$s_{x,i}\colon U_{x,i}\rightarrow G_x$ has image in $W$. Then the monodromy groupoid $M(G,W)$ has a
topology turning it into a topological groupoid.
\end{Cor}

 No now give a result which is similar to  Theorem \ref{comploctr} under  strong conditions.
\begin{Theo}\label{topong} Let X be a topological space and  G  groupoid on $X$. Let W be a compatible
locally trivial subgroupoid of $G$.  Then $G$ has a topological groupoid making it a topological 
groupoid such that $W$ is open in $G$.
\end{Theo}
\begin{Prf} By Theorem \ref{comploctr}, $G$ has a topology  making it a topological groupoid. 
Let $a\in W$ with $\alpha(a)=x$ and $\beta(b)=y$. Since for $x\in O_G$, $s_{x,i}(x)=1_x$, 
$(\wtilde{U}_{x,i})^{-1}a(\wtilde{U}_{y,j})$ is an open neighbourhood of $a$ and since 
$W$ is a subgroupoid,  
$(\wtilde{U}_{x,i})^{-1}a(\wtilde{U}_{y,j})\subseteq W$.  So $W$ is open in $G$.
\end{Prf}

\begin{Cor}
Let G be a  topological groupoid and let  $W$ be a compatible locally trivial subgroupoid of $G$.  Then the monodromy groupoid 
$M(G,W)$ has a topology turning into a topological groupoid such that $\wtilde{\imath}(W)=\wtilde W$ is open
in $M(G,W)$ and any  continuous   pregroupoid morphism  $f\colon W\rightarrow H$  globalises to a
morphism of topological groupoids $\wtilde f\colon M(G,W)\rightarrow H$. 
\end{Cor}
\begin{Prf} Since $\wtilde{\imath}\colon W\rightarrow  \wtilde W$ is a homeomorphism and
$O_G=O_{M(G,W)}$, $\wtilde W$ is a compatible locally trivial subgroupoid of  $M(G,W)$. 
So by Theorem \ref{topong}, $M(G,W)$ becomes a topological groupoid such that $\wtilde W$ 
is open $M(G,W)$. Let $f\colon W\rightarrow H$ be a continuous pregroupoid morphism. By the universal 
property of $M(G,W)$, $f$ globalises to a morphism $\wtilde f\colon M(G,W)\rightarrow H$ which is
continuous on $\wtilde W$. But $\wtilde W$ is open in $M(G,W)$ and $\wtilde W$  generates $M(G,W)$. So the 
 continuity of $\wtilde f$ on whole $M(G,W)$ follows.
\end{Prf}
  
\begin{Theo}\label{top}
Let G be a  topological groupoid which is compatible locally trivial such 
that each star $G_x$ has a universal covering  and let  $W$ be an open subset of
 $G$ such that $O_G\subseteq W$ and $W^2$ is contained in an open 
 neighbourhood $V$ of $O_G$ such that each star $V_x$ is liftable. Suppose that each local
 section $s_{x,i}\colon U_{x,i}\rightarrow G_x$ has image in $W$. Then the monodromy groupoid 
 $M(G,W)$ may be given a topological groupoid structure such that each star 
  $M(G,W)_x$ is isomorphic to the universal covering of $G_x$.
\end{Theo}
\begin{Prf} By Proposition \ref{loctrmon}, the groupoid $M(G,W)$ is compatible locally
 trivial and by Theorem \ref{comploctr} it becomes a topological groupoid. Further by Theorem 
 \ref{startop} we identify $M(G,W)$ with $\Pi G$,  the union of all the universal coverings of the
  stars ${G_x}'s$. 
\end{Prf}

As a corollary of these results we give the following.
\begin{Theo}\label{strmonoth} {\em (monodromy principle)}
Let  $G$ be a  topological groupoid in which each  star $G_x$ is path connected and has a
universal covering   and let $W$ be a locally compatible
subgroupoid of $O_G$ such that $W$ is star liftable. 
 Then any continuous pregroupoid morphism  $f: W \rightarrow H$  determines uniquely a morphism
 $f'\colon \Pi G \rightarrow H$ of topological groupoids such that 
 $f'\wtilde{j} = f$.  
 \hfill $\Box$ 
\end{Theo}

{\bf Acknowledgement}: We would like to thank Prof. Ronald Brown for introducing
this area to us and his help and encouragement.

\end{document}